# Enumeration of standard Young tableaux of shifted strips with constant width

Ping Sun*

*Department of Mathematics, Northeastern University, Shenyang, 110004, China*

**Abstract** Let $g_{n_1,n_2}$ be the number of standard Young tableau of truncated shifted shape with $n_1$ rows and $n_2$ boxes in each row. By using of the integral method this paper derives the recurrence relations of $g_{3,n}$, $g_{n,4}$ and $g_{n,5}$ respectively. Specially, $g_{n,4}$ is the $(2n-1)$-st Pell number.

*MSC*: 05A15, 05E15

*Keywords*: Standard Young tableau; multiple integrals; recurrences; Pell numbers

## 1   Introduction

A shifted diagram of shape $\lambda = (\lambda_1, \cdots, \lambda_d)$ $(\lambda_1 > \cdots > \lambda_d)$ is an array of $|\lambda|$ boxes, where row $i$ (from top to bottom) containing $\lambda_i$ boxes starts with its leftmost box in position $(i,i)$. A standard shifted Young tableau of shape $\lambda$ is a labeling by $\{1, 2, \cdots, |\lambda|\}$ of the boxes in the shifted diagram such that each row and column is increasing (from left to right and from top to bottom respectively). Specially, a standard Young tableaux (SYT) of shifted staircase shape is $\delta_n = (n, n-1, \cdots, 1)$. The enumeration of SYT is an important problem in enumerative combinatorics. See R. P. Stanley's monograph [9] and R. M. Adin and Y. Roichman's recent survey paper [2].

The number of SYT of shifted shape $\lambda$ is given by the well-known product formula[7, 11]:

$$g^\lambda = \frac{|\lambda|!}{\prod_{i=1}^d \lambda_i!} \prod_{i<j} \frac{\lambda_i - \lambda_j}{\lambda_i + \lambda_j}. \qquad (1)$$

---

*E-mail address: plsun@mail.neu.edu.cn

A SYT of truncated shape is the SYT with some boxes removed from the NE corner, which was recently discussed by some authors [1, 6, 10]. Adin et al. in [1] derived the formulas of $\delta_n$ truncated by a square or nearly a square by the method of pivoting theory. G. Panova independently obtained the product formulas of $\delta_n$ truncated by a box in term of Schur function [6].

This paper considers the SYT of shifted shape in case of $\lambda_1 = \cdots = \lambda_d$, which is the SYT of shifted shape $(n_1 + n_2 - 1, n_1 + n_2 - 2, \cdots, n_2)$ truncated by a staircase $\delta_{n_1-1}$, namely a SYT of truncated shifted shape with $n_1$ rows and $n_2$ boxes in each row, illustrated as follows

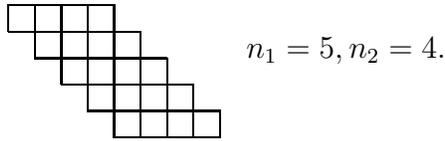
$n_1 = 5, n_2 = 4.$

Let $g_{n_1,n_2}(n_1, n_2 \geq 2)$ be the number of SYT of truncated shifted shape with $n_1$ rows and $n_2$ boxes in each row. It is clear that $g_{2,n}$ is the $(n-1)$-th Catalan number $\frac{1}{n}\binom{2n-2}{n-1}$, $g_{n,2} = 1$ and $g_{n,3} = 2^{n-1}$. This kind of SYT of shifted strips with constant width has many applications. J. B. Lewis conjectured that the number of alternating permutations of length $2n - 2$ avoiding the pattern 3412 is $g_{3,n}$, which summation representation was given by G. Pabova [5]. In fact, $g_{n_1,n_2}$ is also the number of $n_1 \times n_2$ matrices containing a permutation of $[n_1n_2]$ in increasing order rowwise, columnwise, diagonally and (downwards) antidiagonally, because the character in increasing order antidiagonally of matrix coincides with the shifted property of SYT.

Suppose $T(n, k)$ is the number of $n \times k$ matrices containing a permutation of $[nk]$ in increasing order rowwise, columnwise, diagonally and antidiagonally, R. H. Hardin gives several empirical formulas of $T(n, k)$ in case of $k \leq 7$ [8, A181196]:

**Empirical Recurrence Relations** (*R.H.Hardin*)

$T(n, 4): a_n = 6a_{n-1} - a_{n-2}.$ (2)

$T(n, 5): a_n = 24a_{n-1} - 40a_{n-2} - 8a_{n-3}.$ (3)

$T(n, 6): a_n = 120a_{n-1} - 1672a_{n-2} + 544a_{n-3} - 6672a_{n-4} + 256a_{n-5}.$ (4)

$T(n, 7): a_n = 720a_{n-1} - 84448a_{n-2} + 1503360a_{n-3} - 17912224a_{n-4} - 318223104a_{n-5}$
$\qquad + 564996096a_{n-6} + 270471168a_{n-7} - 11373824a_{n-8} + 65536a_{n-9}.$ (5)

In addition, Ping Sun in [10] shows that $g_{n_1,n_2}$ is involved in the nested distribution of $n_1$ groups of independent order statistics with $n_2$ samples from uniform distribution on interval $(0,1)$. Generally, the order statistics model of SYT in [10] implies

**Proposition 1**. *For $n_1, n_2 \geq 2$,*

$$g_{n_1,n_2} = (n_1 n_2)! \cdot J(n_1, n_2) = (n_1 n_2)! \int \cdots \int_{D(n_1,n_2)} dx_{i,j}, \tag{6}$$

*where $D(n_1, n_2)$ is the following SYT-type integral domain (the variables are increasing from left to right and from top to bottom)*

$$
\begin{array}{cccc}
0 < x_{1,1} < & x_{1,2} & < \cdots & < x_{1,n_2} \\
& \wedge & & \wedge \\
& x_{2,1} & < \cdots < \quad x_{2,n_2-1} & < x_{2,n_2} \\
& & \ddots \quad \ddots & \ddots \\
& & x_{n_1,1} < \cdots \quad \cdots & < x_{n_1,n_2} < 1.
\end{array}
$$

We derive the recurrence formulas of $g_{3,n}, g_{n,4}$ and $g_{n,5}$ by using of the integral method of [10] in this paper. In Section 2 we evaluate the nested distribution of *three* groups of independent order statistics with $n$ samples from uniform distribution on interval $(0, 1)$, which implies a new summation representation of $g_{3,n}$. So that a non-homogeneous linear recurrence relation of $g_{3,n}$ is given. In Section 3 we compute the corresponding multiple integrals and prove the recurrence relations (2) and (3) respectively. In particular, $g_{n,4}$ is shown to be the $(2n-1)$-st Pell number which implies the empirical formula (2).

## 2 New recurrence relation of $g_{3,n}$

There are two results of $g_{3,n}$ in the literature. Considering the enumeration of SYT of shifted shape $(n, n-1, i), 0 \leq i \leq n-2$, G. Panova gives the following summation representation.

**Proposition 2**. [5] *For $n \geq 2$,*

$$g_{3,n} = \sum_{i=0}^{n-2} \frac{(2n+i-1)!(n-i)(n-i-1)}{n! \, (n-1)! \, i! \, (2n-1)(n+i)(n+i-1)}. \tag{7}$$

For the number of $3 \times n$ matrices containing a permutation of $[3n]$ in increasing order rowwise, columnwise, diagonally and antidiagonally, V. Kotesovec gives the complicated order 2 recurrence relation.

**Proposition 3.** [8, A181197] *For $n \geq 3$,*

$$(2n-1)(7n-13)n^2 \cdot g_{3,n} = 2(182n^4 - 1185n^3 + 2722n^2 - 2625n + 900) \cdot g_{3,n-1}$$
$$+ 3(2n-5)(3n-5)(3n-4)(7n-6) \cdot g_{3,n-2}, \quad g_{3,1} = 1, \ g_{3,2} = 1. \quad (8)$$

For $0 < t_1 < t_2 < t_3 < 1$, we consider the following integral

$$J_{n-1}(t_1, t_2, t_3) = \int \cdots \int_{D(t)} dx_{i,j}, n \geq 2,$$

where $D(t)$ is the following SYT-type integral domain

$$0 < x_{i,1} < x_{i,2} < \cdots < x_{i,n-1} < t_i, 1 \leq i \leq 3; \quad x_{1,2} < x_{2,1}, x_{2,n-1} < x_{3,n-2};$$
$$t_1 < x_{2,n-1}, t_2 < x_{3,n-1}; \quad x_{1,j+2} < x_{2,j+1} < x_{3,j}, 1 \leq j \leq n-3.$$

From Proposition 1, there is

$$g_{3,n} = (3n)! \iiint\limits_{0<t_1<t_2<t_3<1} J_{n-1}(t_1, t_2, t_3) dt_i, \quad n \geq 2. \quad (9)$$

**Lemma 1.** *For $n \geq 1$, $0 < t_1 < t_2 < t_3 < 1$,*

$$J_n(t_1, t_2, t_3) = \frac{t_1^n t_2^n t_3^n}{n!^3} + \frac{2}{n!(2n)!} \sum_{i=0}^{n-1} (-1)^{n-i} \binom{2n}{i} \left[ t_1^n t_2^{2n-i} t_3^i - t_1^{2n-i} t_2^n t_3^i + t_1^{2n-i} t_2^i t_3^n \right]. (10)$$

**Proof.** It is clear that

$$J_1(t_1, t_2, t_3) = t_1(t_2 - t_1)(t_3 - t_2) = t_1 t_2 t_3 - (t_1 t_2^2 - t_1^2 t_2 + t_1^2 t_3).$$

Suppose (10) is true for $n-1$, then

$$J_n(t_1, t_2, t_3) = \iiint\limits_{0<x<t_1<y<t_2<z<t_3} J_{n-1}(x, y, z) dx dy dz$$

$$= \frac{t_1^n(t_2^n - t_1^n)(t_3^n - t_2^n)}{n!^3} + \frac{2}{(n-1)!(2n-2)!} \sum_{i=0}^{n-2} (-1)^{n-1-i} \binom{2n-2}{i} \frac{1}{n(i+1)(2n-i-1)} \times$$

$$\left[ t_1^n (t_2^{2n-i-1} - t_1^{2n-i-1})(t_3^{i+1} - t_2^{i+1}) - t_1^{2n-i-1}(t_2^n - t_1^n)(t_3^{i+1} - t_2^{i+1}) + t_1^{2n-i-1}(t_2^{i+1} - t_1^{i+1})(t_3^n - t_2^n) \right]$$

$$= \frac{t_1^n t_2^n t_3^n}{n!^3} + \frac{2}{n!(2n)!} \sum_{i=0}^{n-1} (-1)^{n-i} \binom{2n}{i} \left[ t_1^n t_2^{2n-i} t_3^i - t_1^{2n-i} t_2^n t_3^i + t_1^{2n-i} t_2^i t_3^n \right] + R_n(t_1, t_2, t_3),$$

where

$$R_n(t_1, t_2, t_3) = (t_1^{2n}t_2^n - t_1^n t_2^{2n} - t_1^{2n}t_3^n)\left[\frac{1}{n!^3} + \frac{2}{n!(2n)!}\sum_{i=0}^{n-1}(-1)^{n-i}\binom{2n}{i}\right] = 0$$

follows from the known identity [4, 1.86]

$$\sum_{i=0}^{n}(-1)^i\binom{2n}{i} = (-1)^n\binom{2n-1}{n}. \tag{11}$$

The proof of lemma 1 is complete by induction. □

It should be noted that $J_{n-1}(t_1, t_2, t_3)$ is the shifted nested conditional distribution of *three* groups of independent order statistics with $n$ samples from uniform distribution on interval $(0, 1)$. The following result gives a non-homogeneous linear recurrence of $g_{3,n}$.

**Theorem 1.** *For $n \geq 1$, the number $g_{3,n}$ of SYT of truncated shifted shape with 3 rows and $n$ boxes in each row satisfies*

$$g_{3,n+1} = -g_{3,n} + \frac{7n+1}{n^2(n+1)^2}\binom{2n-2}{n-1}\binom{3n}{n-1}, \quad g_{3,1} = 1. \tag{12}$$

**Proof.** For $n \geq 2$, combining (9) and (10) we have the summation representation

$$g_{3,n} = \frac{(3n)!}{6 \cdot n!^3} + \sum_{i=0}^{n-2}(-1)^{n-1-i}\frac{(3n-1)!(5n-3i-3)}{n!\, i!\, (2n-i-1)!\, (3n-i-1)}. \tag{13}$$

Decomposing $5n - 3i - 3 = 3(3n - i - 1) - 4n$ in above summation, from (11) and Frisch's identity [4, 4.2]:

$$\sum_{i=0}^{2n-1}(-1)^i\binom{2n-1}{i}\frac{1}{n+i} = \frac{(n-1)!(2n-1)!}{(3n-1)!},$$

the summation representation (13) of $g_{3,n}$ is simplified to be

$$g_{3,n} = -\frac{4n-5}{6(2n-1)} \cdot \frac{(3n)!}{n!^3} + 4(-1)^{n-1} + (-1)^n 4n\binom{3n-1}{n} \cdot A_n, \tag{14}$$

where

$$A_n = \sum_{i=0}^{n}(-1)^i\binom{2n-1}{i}\frac{1}{n+i}.$$

The recurrence of $A_n$ is not difficult to derive by using of the equality (11).

$$A_{n+1} = \frac{1}{3n+2}\sum_{i=0}^{n+1}(-1)^i\frac{(2n+1)!}{i!\,(2n-i)!}\left[\frac{1}{2n+1-i}+\frac{1}{n+1+i}\right]$$

$$= \frac{(-1)^{n+1}\binom{2n}{n+1}}{3n+2} + \frac{2n(2n+1)}{(3n+1)(3n+2)}\sum_{i=0}^{n+1}(-1)^i\binom{2n-1}{i}\left[\frac{1}{2n-i}+\frac{1}{n+1+i}\right]$$

$$= \frac{(-1)^{n+1}(10n^3+8n^2-n-1)\binom{2n}{n+1}}{2n(n+1)(3n+1)(3n+2)} + \frac{2n(2n+1)}{(3n+1)(3n+2)}\sum_{i=0}^{n}(-1)^i\frac{\binom{2n-1}{i}}{n+1+i},$$

and the summation in last equality is equal to

$$\sum_{i=0}^{n}(-1)^{i+1}\frac{1}{n+1+i}\left[\binom{2n-1}{i+1}-\binom{2n}{i+1}\right]$$

$$= \frac{(-1)^{n+1}\left[\binom{2n-1}{n+1}-\binom{2n}{n+1}\right]}{2n+1} + A_n - \frac{1}{3n}\sum_{i=0}^{n}(-1)^i\frac{(2n)!}{i!\,(2n-1-i)!}\left[\frac{1}{n+i}+\frac{1}{2n-i}\right]$$

$$= \frac{1}{3}A_n - \frac{n^2-1}{6n^2(2n+1)}(-1)^{n+1}\binom{2n}{n+1},$$

therefore the recurrence of $A_n$ is

$$A_{n+1} = \frac{2n(2n+1)}{3(3n+1)(3n+2)}A_n + (-1)^{n+1}\frac{28n^3+22n^2-n-1}{6n(n+1)(3n+1)(3n+2)}\binom{2n}{n+1},\quad A_1 = \frac{1}{2}.$$

So that the recurrence relation (12) of $g_{3,n}$ follows from

$$A_n = \frac{(-1)^n}{4n\binom{3n-1}{n}}g_{3,n} + \frac{1}{n\binom{3n-1}{n}} + (-1)^n\frac{4n-5}{8n(2n-1)}\binom{2n-1}{n}.$$

The proof of theorem 1 is complete. □

## 3 Recurrence relations of $g_{n,4}$ and $g_{n,5}$

It is well-known the Pell numbers $P_n$ are defined by the recurrence relation [8, A000129]

$$P_n = 2P_{n-1} + P_{n-2},\quad P_0 = 0, P_1 = 1. \tag{15}$$

The Pell numbers arise historically in the rational approximation to $\sqrt{2}$, which are used to enumerate the numbers of certain pattern-avoiding permutations recently[3].

**Theorem 2**. *For $n \geq 1$, the number $g_{n,4}$ of SYT of truncated shifted shape with $n$ rows and 4 boxes in each row is the $(2n-1)$-st Pell number $P_{2n-1}$, which satisfies*

$$g_{n,4} = 6g_{n-1,4} - g_{n-2,4}, \quad g_{1,4} = 1, g_{2,4} = 5. \tag{16}$$

**Proof.** For $1 \leq i \leq n$, write the variable $t_{2i-1}$ corresponding to the box $(i, i+1)$, $t_{2i}$ corresponding to the box $(i, i+2)$ in the SYT of shifted strip with width 4 respectively:

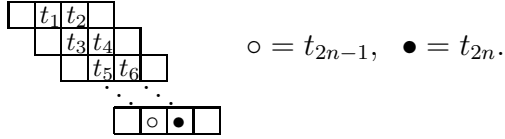
$\circ = t_{2n-1}, \quad \bullet = t_{2n}.$

Proposition 1 implies that

$$g_{n,4} = (4n)! \int \cdots \int_{0<t_1<t_2<t_3<\cdots<t_{2n}<1} t_1(1-t_{2n}) \prod_{i=2}^{n}(t_{2i-1}-t_{2i-3})(t_{2i}-t_{2i-2})dt_1\cdots t_{2n}$$

$$= (4n)! \iint_{0<t_{2n-1}<t_{2n}<1} (1-t_{2n})J_n(t_{2n-1},t_{2n})dt_{2n-1}dt_{2n}, \quad n \geq 2.$$

We shall now use the method of induction to prove the following

$$J_n(t_{2n-1}, t_{2n}) = \int \cdots \int_{0<t_1<t_2<t_3<\cdots<t_{2n-2}<t_{2n-1}} t_1 \prod_{i=2}^{n}(t_{2i-1}-t_{2i-3})(t_{2i}-t_{2i-2})dt_1\cdots t_{2n-2}$$

$$= \frac{t_{2n-1}^{4n-4}}{(4n-3)!}\left\{(4n-3)P_{2n-2}t_{2n} - [(4n-4)P_{2n-2} - P_{2n-3}]t_{2n-1}\right\}, \tag{17}$$

where $P_i$ is the Pell number.

It is clear that

$$J_2(t_3, t_4) = \iint_{0<t_1<t_2<t_3} t_1(t_3-t_1)(t_4-t_2)dt_1dt_2 = \frac{t_3^4}{5!}(10t_4 - 7t_3),$$

which agrees with (17) because $P_1 = 1, P_2 = 2$.

Furthermore, from the recurrence relation (15) of Pell numbers,

$$J_{n+1}(t_{2n+1}, t_{2n+2}) = \iint_{0<t_{2n-1}<t_{2n}<t_{2n+1}} J_n(t_{2n-1}, t_{2n})(t_{2n+1} - t_{2n-1})(t_{2n+2} - t_{2n})dt_{2n-1}dt_{2n}$$

$$= \int_0^{t_{2n+1}} \frac{t_{2n}^{4n-2}}{(4n-1)!}\{(4n-1)P_{2n-1}t_{2n+1} - [(4n-2)P_{2n-1} - P_{2n-2}]t_{2n}\}(t_{2n+2} - t_{2n})dt_{2n}$$

$$= \frac{t_{2n+1}^{4n}}{(4n+1)!}[(4n+1)P_{2n}t_{2n+2} - (4nP_{2n} - P_{2n-1})t_{2n+1}],$$

which shows (17) is true. Therefore,

$$g_{n,4} = (4n)! \iint_{0<t_{2n-1}<t_{2n}<1} (1 - t_{2n})J_n(t_{2n-1}, t_{2n})dt_{2n-1}dt_{2n}$$

$$= \frac{(4n)!}{(4n-3)!} \int_0^1 \left[P_{2n-2} - \frac{(4n-4)P_{2n-2} - P_{2n-3}}{4n-2}\right] t_{2n}^{4n-2}(1 - t_{2n})dt_{2n}$$

$$= 4n(4n-1) \int_0^1 P_{2n-1}t_{2n}^{4n-2}(1 - t_{2n})dt_{2n} = P_{2n-1}.$$

It is clear (16) follows from the recurrence relation of Pell number. $\square$

**Theorem 3**. *For $n \geq 4$, the numbers $g_{n,5}$ of SYT of truncated shifted shape with $n$ rows and 5 boxes in each row satisfy the following recurrence relation*

$$g_{n,5} = 24g_{n-1,5} - 40g_{n-2,5} - 8g_{n-3,5}, \quad g_{1,5} = 1, g_{2,5} = 14, g_{3,5} = 290. \tag{18}$$

**Proof.** We shall derive the recurrence of $g_{n,5}$ from the relations of certain integrals. For convenient, write the variables $0 < x_i < y_i < z_i < s_i < t_i < 1$ corresponding to the five boxes in row $i$ ($1 \leq i \leq n$) in the SYT of shifted strip with width 5 respectively.

$$\begin{array}{ccccc}
\ddots & \ddots & \ddots & \ddots & \ddots \\
\boxed{x_{n-1}\ y_{n-1}\ z_{n-1}\ s_{n-1}\ t_{n-1}} & & & & \\
& \boxed{x_n\ \ y_n\ \ z_n\ \ s_n\ \ t_n} & & &
\end{array}$$

From Proposition 1, we have

$$g_{n,5} = (5n)! \int \cdots \int_{D(n,5)} dx_i dy_i dz_i ds_i dt_i = (5n)! \cdot J_n(5).$$

Denote $D_1(x_1, y_1, z_1, s_1) = 1$, consider the following integral

$$D_n(x_n, y_n, z_n, s_n) = \int \cdots \int_{\substack{D_{n-1,5},\, y_{n-1}<x_n,\\ z_{n-1}<y_n,\, s_{n-1}<z_n,\, t_{n-1}<s_n}} dx_1 y_1 z_1 s_1 t_1 \cdots dx_{n-1} y_{n-1} z_{n-1} s_{n-1} t_{n-1},$$

$J_n(5)$ can be written to be

$$J_n(5) = \int \cdots \int_{\substack{0<x_{n-1}<y_{n-1}<z_{n-1}<s_{n-1}<t_{n-1}<s_n,\\ s_{n-1}<z_n,\, z_{n-1}<y_n<z_n<s_n<t_n<1}} (y_n - y_{n-1}) D_{n-1}(x_{n-1}, \cdots, s_{n-1}) dx_{n-1} \cdots t_{n-1} dy_n z_n s_n t_n$$

$$= \int \cdots \int_{\substack{0<z_{n-1}<s_{n-1}<t_{n-1}<s_n,\\ s_{n-1}<z_n,\, z_{n-1}<y_n<z_n<s_n<t_n<1}} (A_n y_n - B_n) dz_{n-1} s_{n-1} t_{n-1} dy_n z_n s_n t_n, \tag{19}$$

where

$$A_n = \iint_{0<x_{n-1}<y_{n-1}<z_{n-1}} D_{n-1}(x_{n-1}, y_{n-1}, z_{n-1}, s_{n-1}) dx_{n-1} dy_{n-1}$$

$$= C_1(n) \frac{z_{n-1}^{5n-9}}{(5n-9)!} s_{n-1} - C_2(n) \frac{z_{n-1}^{5n-8}}{(5n-8)!},$$

$$B_n = \iint_{0<x_{n-1}<y_{n-1}<z_{n-1}} y_{n-1} D_{n-1}(x_{n-1}, y_{n-1}, z_{n-1}, s_{n-1}) dx_{n-1} dy_{n-1}$$

$$= C_3(n) \frac{z_{n-1}^{5n-8}}{(5n-8)!} s_{n-1} - C_4(n) \frac{z_{n-1}^{5n-7}}{(5n-7)!}.$$

Notice that the definition of $D_n(x_n, y_n, z_n, s_n)$ implies

$$A_{n+1} = \iint_{0<x_n<y_n<z_n} D_n(x_n, y_n, z_n, s_n) dx_n dy_n$$

$$= \int \cdots \int_{\substack{0<z_{n-1}<s_{n-1}<t_{n-1}<s_n,\\ s_{n-1}<z_n,\, z_{n-1}<y_n<z_n<s_n}} (A_n y_n - B_n) dz_{n-1} s_{n-1} t_{n-1} dy_n, \tag{20}$$

then

$$A_{n+1} = [10(n-1)(5n-7)C_1(n) - (10n-11)C_2(n) - (10n-11)C_3(n) + 2C_4(n)] \times$$

$$\frac{z_n^{5n-4}}{(5n-4)!} s_n - [(5n-4)(5n-7)(10n-11)C_1(n) - 2(5n-4)(5n-6)C_2(n)$$

$$- 50(n-1)^2 C_3(n) + (10n-9)C_4(n)] \frac{z_n^{5n-3}}{(5n-3)!}.$$

By the similar arguments,

$$B_{n+1} = \iint\limits_{0<x_n<y_n<z_n} y_n D_n(x_n, \cdots, s_n) dx_n dy_n$$

$$= \int \cdots \int\limits_{\substack{0<z_{n-1}<s_{n-1}<t_{n-1}<s_n, \\ s_{n-1}<z_n, z_{n-1}<y_n<z_n<s_n}} (A_n y_n^2 - B_n y_n) dz_{n-1} s_{n-1} t_{n-1} y_n$$

$$= [(5n-4)(5n-7)(10n-11)C_1(n) - 50(n-1)^2 C_2(n) - 2(5n-4)(5n-6)C_3(n)$$
$$+ (10n-9)C_4(n)] \frac{z_n^{5n-3}}{(5n-3)!} s_n - (5n-3)[2(5n-4)(5n-6)(5n-7)C_1(n)$$
$$- (5n-6)(10n-9)C_2(n) - (5n-6)(10n-9)C_3(n) + 10(n-1)C_4(n)] \frac{z_n^{5n-2}}{(5n-2)!}.$$

Therefore, For $n \geq 2$, the recurrence relations of $C_i(n)$ $(1 \leq i \leq 4)$ are

$$C_1(n+1) = 10(n-1)(5n-7)C_1(n) - (10n-11)C_2(n) - (10n-11)C_3(n) + 2C_4(n),$$
$$C_2(n+1) = (5n-4)(5n-7)(10n-11)C_1(n) - 2(5n-4)(5n-6)C_2(n)$$
$$\qquad - 50(n-1)^2 C_3(n) + (10n-9)C_4(n),$$
$$C_3(n+1) = (5n-4)(5n-7)(10n-11)C_1(n) - 50(n-1)^2 C_2(n)$$
$$\qquad - 2(5n-4)(5n-6)C_3(n) + (10n-9)C_4(n),$$
$$C_4(n+1) = 2(5n-3)(5n-4)(5n-6)(5n-7)C_1(n) - (5n-3)(5n-6)(10n-9)C_2(n)$$
$$\qquad - (5n-3)(5n-6)(10n-9)C_3(n) + 10(n-1)(5n-3)C_4(n),$$

with the initial values $C_1(2) = 0$, $C_2(2) = -1$, $C_3(2) = 0$, $C_4(2) = -2$.

On the other hand, combining (19) and (20), we have

$$J_n(5) = \iiint\limits_{0<z_n<s_n<t_n<1} A_{n+1} dz_n s_n t_n$$

$$= \iiint\limits_{0<z_n<s_n<t_n<1} [C_1(n+1) \frac{z_n^{5n-4}}{(5n-4)!} s_n - C_2(n+1) \frac{z_n^{5n-3}}{(5n-3)!}] dz_n s_n t_n$$

$$= \frac{(5n-2)C_1(n+1) - C_2(n+1)}{(5n)!},$$

then,

$$g_{n,5} = (5n-2)C_1(n+1) - C_2(n+1), \tag{21}$$

and

$$g_{n,5} = (5n-7)(25n-24)C_1(n) - (25n-26)C_2(n) - (25n-28)C_3(n) + 5C_4(n), \quad (22)$$

which follows from (21) and the recurrences of $C_i(n)$.

Furthermore, by using of the recurrence relations of $C_i(n)$ again, (22) implies

$$\begin{aligned}g_{n+1,5} =&(5n-2)(25n+1)C_1(n+1) - (25n-1)C_2(n+1) \\ &- (25n-3)C_3(n+1) + 5C_4(n+1) \\ =&(5n-7)(550n-524)C_1(n) - (550n-590)C_2(n) \\ &- (550n-594)C_3(n) + 110C_4(n).\end{aligned} \quad (23)$$

So that from (21)-(23) we have

$$\begin{aligned}g_{n+1,5} &= 22g_{n,5} + 4(5n-7)C_1(n) + 18C_2(n) - 22C_3(n) \\ &= 22g_{n,5} + 4g_{n-1,5} + 22[C_2(n) - C_3(n)], \quad n \geq 2.\end{aligned}$$

Finally, the recurrence relation of $C_i(n)$ shows

$$C_2(n+1) - C_3(n+1) = 2[C_2(n) - C_3(n)], \quad n \geq 2,$$

which implies

$$g_{n+1,5} - 22g_{n,5} - 4g_{n-1,5} = 2(g_{n,5} - 22g_{n-1,5} - 4g_{n-2,5}), \quad n \geq 3.$$

The initial values $g_{1,5} = 1, g_{2,5} = 14$ and $g_{3,5} = 290$ are obvious, then the proof of theorem 3 is complete. □

# Acknowledgements

The author would like to thank Yu Bao and Yanting He for helpful discussions.

# References

[1] Ron M. Adin, Ronald C. King, and Yuval Roichman, Enumeration of standard Young tableaux of certain truncated shapes. Electron. J. Combin., 18(2) (2011), # P20.

[2] Ron M. Adin, Yuval Roichman, Enumeration of Standard Young Tableaux. arXiv:1408.4497v2, 2014.

[3] M. Barnabei, F. Bonetti, and M. Silimbani, Two permutation classes related to the Bubble Sort operator. Electron. J. Combin., 19(3) (2012), # P25.

[4] H. W. Gould, Combinatorial Identities. Morgantown, W. Va. 1972.

[5] J. B. Lewis, Generating trees and pattern avoidance in alternating permutations. Electron. J. Combin., 19 (2012), # P21.

[6] Greta Panova, Tableaux and plane partitions of truncated shapes. Adv. Appl. Math., 49:196–217, 2012.

[7] I. Schur, On the representation of the symmetric and alternating groups by fractional linear substitutions. J. Reine Angew. Math. 139:155–250, 1911.

[8] N. J. A. Sloane, The On-Line Encyclopedia of Integer Sequences. http://oeis.org.

[9] Richard P. Stanley, Enumerative Combinatorics. vol. 2, Cambridge University Press, New York, 1999.

[10] Ping Sun, Evaluating the numbers of some skew standard Young tableaux of truncated shapes. Electron. J. Combin. 22(1) (2015), # P1.2.

[11] R. M. Thrall, A combinatorial problem. Michigan Math. J. 1:81–88, 1952.